\documentclass{amsart}
\usepackage[leqno]{amsmath}
\usepackage{amssymb}
\usepackage{amsthm}
\usepackage{amsfonts}
\usepackage{mathrsfs}
\usepackage{mathalfa}
\usepackage{enumerate}
\usepackage{xcolor}
\usepackage{subfigure}
\usepackage{morefloats}
\usepackage{graphicx}
\usepackage{hyperref}
\usepackage{caption}
\usepackage{multicol}
\usepackage{enumitem}

\numberwithin{equation}{section}

\newtheorem{theorem}{Theorem}[section]

\newtheorem{lemma}[theorem]{Lemma}

\newcommand{\C}{\mathbb{C}}
\newcommand{\D}{\mathbb{D}}

\newcommand{\h}{\mathbb{H}}
\newcommand{\N}{\mathbb{N}}
\newcommand{\Z}{\mathbb{Z}}

\newcommand{\R}{\mathbb{R}}
\newcommand{\s}{\mathbb{S}}

\newcommand{\cD}{\mathcal {D}}

\newcommand{\cQ}{\mathcal {Q}}

\newcommand{\cS}{\mathcal {S}}
\newcommand{\cT}{\mathcal {T}}

\newcommand{\cW}{\mathcal {W}}

\newcommand{\SLE}{{\rm SLE}}
\newcommand{\QLE}{{\rm QLE}}
\newcommand{\CLE}{{\rm CLE}}

\newcommand{\re}{\text{Re}}
\newcommand{\cov}{\text{cov}}

\newcommand{\wt}{\widetilde}
\newcommand{\wh}{\widehat}
\newcommand{\ol}{\overline}

\newcommand{\sol}[1]{{}}

\newcommand{\strip}{{\mathscr{S}}}
\newcommand{\cyl}{\mathscr{C}}

\newcommand{\mustwo}{\mu_{{\mathrm{SPH}}}^2}
\newcommand{\fb}[2]{B^\bullet(#1,#2)}

\begin{document}

\title{Liouville quantum gravity as~a~metric~space~and~a~scaling~limit}
\author{Jason Miller}

\date{\today}
\maketitle

\vspace{-0.05\textheight}

\begin{abstract}
Over the past few decades, two natural random surface models have emerged within physics and mathematics. The first is Liouville quantum gravity, which has its roots in string theory and conformal field theory from the 1980s and 1990s. The second is the Brownian map, which has its roots in planar map combinatorics from the 1960s together with recent scaling limit results.  This article surveys a series of works with Sheffield in which it is shown that Liouville quantum gravity (LQG) with parameter $\gamma=\sqrt{8/3}$ is equivalent to the Brownian map.  We also briefly describe a series of works with Gwynne which use the $\sqrt{8/3}$-LQG metric to prove the convergence of self-avoiding walks and percolation on random planar maps towards $\SLE_{8/3}$ and $\SLE_6$, respectively, on a Brownian surface.
\end{abstract}

\section{Introduction}

\subsection{The Gaussian free field}

Suppose that $D \subseteq \C$ is a planar domain.  Informally, the Gaussian free field (GFF) $h$ on $D$ is a Gaussian random variable with covariance function
\[ \cov(h(x),h(y)) = G(x,y) \quad\text{for}\quad x,y \in D\]
where~$G$ denotes the Green's function for~$\Delta$ on~$D$.  Since $G(x,y) \sim -\log|x-y|$ as $x \to y$, it turns out that it is not possible to make sense of the GFF as a function on~$D$ but rather it exists as a distribution on~$D$.  Perhaps the most natural construction of the GFF is using a series expansion.  More precisely, one defines $H_0^1(D)$ to be the Hilbert space closure of $C_0^\infty(D)$ with respect to the \emph{Dirichlet inner product}
\begin{equation}
\label{eqn:dirichlet}
	(f,g)_\nabla = \frac{1}{2\pi} \int_D \nabla f(x) \cdot \nabla g(x) dx.
\end{equation}
One then sets
\begin{equation}
\label{eqn:gff_series}
	h = \sum_n \alpha_n \phi_n
\end{equation}
where $(\alpha_n)$ is a sequence of i.i.d.\ $N(0,1)$ random variables and $(\phi_n)$ is an orthonormal basis of $H_0^1(D)$.  The convergence of the series~\eqref{eqn:gff_series} does not take place in $H_0^1(D)$, but rather in the space of distributions on $D$.  (One similarly defines the GFF with free boundary conditions on a given boundary segment $L$ using the same construction but including also the those functions can be non-zero on $L$.)  Since the Dirichlet inner product is conformally invariant, so is the law of the GFF.

The GFF is a fundamental object in probability theory.  Just like the Brownian motion arises as the scaling limit of many different types of random curves, the GFF describes the scaling limit of many different types of random surface models \cite{K,BAD96,gos2001phi,RV08,m2011bounded}.  It also has deep connections with many other important objects in probability theory, such as random walks and Brownian motion \cite{dynkin1,dynkin2,lejan2008markov} and the Schramm-Loewner evolution (SLE) \cite{S0,ss2009contour,SchrammSheffieldGFF2,she2010zipper, dub2009partition,ms2012imag1,ms2012imag2,ms2012imag3,ms2013imag4,matingtrees}.

\subsection{Liouville quantum gravity}

Liouville quantum gravity (LQG) is one of several important random geometries that one can build from the GFF.  It was introduced (non-rigorously) in the physics literature by Polyakov in the 1980s as a model for a string \cite{pol81bosonic,pol81fermionic,pol88qg}.  In its simplest form, it refers to the random Riemannian manifold with metric tensor given by
\begin{equation}
\label{eqn:metric_tensor}
e^{\gamma h(z)}(dx^2 + dy^2)
\end{equation}
where $h$ is (some variant of) the GFF on $D$, $\gamma$ is a real parameter, $z = x+iy \in D$, and $dx^2 +dy^2$ denotes the Euclidean metric on $D$.  This expression does not make literal sense because $h$ is only a distribution on $D$ and not a function, hence does not take values at points.

There has been a considerable effort within the probability community over the course of the last decade or so to make rigorous sense of LQG.  This is motivated by the goal of putting a heuristic from the physics literature, the so-called KPZ formula \cite{kpz1988} which has been used extensively to (non-rigorously) derive critical exponents for two-dimensional lattice models, onto firm mathematical ground.  Another motivation comes from deep conjectures, which we will shortly describe in more detail, which state that LQG should describe the large-scale behavior of random planar maps.

The volume form
\begin{equation}
\label{eqn:volume_form}
e^{\gamma h(z)} dxdy
\end{equation}
associated with the metric~\eqref{eqn:metric_tensor} was constructed in \cite{ds2011kpz} by regularizing $h$ by considering its average $h_\epsilon(z)$ on $\partial B(z,\epsilon)$ and then setting
\begin{equation}
\label{eqn:volume_form_eps}
 \mu_h^\gamma = \lim_{\epsilon \to 0} \epsilon^{\gamma^2/2} e^{\gamma h_\epsilon(z)} dx dy
\end{equation}
where $dx dy$ denotes Lebesgue measure on $D$.  We note that there is no difficulty in making sense of the expression inside of the limit on the right hand side of~\eqref{eqn:volume_form_eps} for each fixed $\epsilon > 0$ since $(z,\epsilon) \mapsto h_\epsilon(z)$ is a continuous function \cite{ds2011kpz}.  The factor $\epsilon^{\gamma^2/2}$ appears because the leading order term in the variance of $h_\epsilon(z)$ is $\log \epsilon^{-1}$.  In the case that $h$ is a GFF on a domain $D$ with a linear boundary segment $L$ and free boundary conditions along $L$, one can similarly construct a boundary length measure by setting
\begin{equation}
\label{eqn:boundary_measure_eps}
\nu_h^\gamma = \lim_{\epsilon \to 0} \epsilon^{\gamma^2/4} e^{\gamma h_\epsilon(x)/2} dx
\end{equation}
where $dx$ denotes Lebesgue measure on $L$.  There is in fact a general theory of random measures with the same law as $\mu_h$, $\nu_h$ which is referred to as \emph{Gaussian multiplicative chaos} and was developed by Kahane \cite{kahane} in the 1980s; see \cite{rhodes-vargas-review} for a more recent review.  (Similar measures also appeared earlier in \cite{hk1971quantum}.)

The regularization procedure used to define the area and boundary measures leads to a natural change of coordinates formula for LQG.  Namely, if $h$ is a GFF on a domain $D$, $\varphi \colon \wt{D} \to D$ is a conformal transformation, and one takes
\begin{equation}
\label{eqn:change_of_coordinates}
\wt{h} = h \circ \varphi + Q \log|\varphi'| \quad\text{where}\quad Q = \frac{2}{\gamma} + \frac{\gamma}{2}
\end{equation}
then the area and boundary measures defined by $\wt{h}$ are the same as the pushforward of the area and boundary measures defined by $h$.  Therefore $(D,h)$ and $(\wt{D},\wt{h})$ can be thought of as parameterizations of the same surface.  Whenever two pairs $(D,h)$ and $(\wt{D},\wt{h})$ are related as in~\eqref{eqn:change_of_coordinates}, we say that they are equivalent as quantum surfaces and a \emph{quantum surface} is an equivalence class under this relation.  A particular choice of representative is referred to as an \emph{embedding}.  There are many different choices of embeddings of a given quantum surface which can be natural depending on the context, and this is a point we will come back to later.  We note that when $\gamma \to 0$ so that an LQG surface corresponds to a flat, Euclidean surface (i.e., the underlying planar domain), the change of coordinate formula~\eqref{eqn:change_of_coordinates} exactly corresponds to the usual change of coordinates formula.

\subsection{Random planar maps}

A \emph{planar map} is a graph together with an embedding into the plane so that no two edges cross.  Planar maps $M_1,M_2$ are said to be equivalent if there exists an orientation preserving homeomorphism $\varphi$ of $\R^2$ which takes $M_1$ to $M_2$.  The \emph{faces} of a planar map are the connected components of the complement of its edges.  A map is a called  a triangulation (resp.\ quadrangulation) if it has the property that all faces have exactly three (resp.\ four) adjacent edges.  The study of planar maps goes back to work of Tutte \cite{tutte1962census}, who worked on enumerating planar maps in the context of the four color theorem, and of Mullin \cite{MR0205882}, who worked on enumerating planar maps decorated with a distinguished spanning tree.  Random maps were intensively studied in physics by random matrix techniques, starting from \cite{bipz1978matrix} (see, e.g., the review \cite{dgz1995review}).  This field has since been revitalized by the development of bijective techniques due to Cori-Vauquelin \cite{cv1981maps} and Schaeffer \cite{schaeffer} and has remained a very active area in combinatorics and probability, especially in the last 20 or so years.

Since there are only a finite number of planar maps with $n$ faces, one can pick one uniformly at random.  We will now mention some of the main recent developments in the study of the metric properties of large uniformly random maps.  First, it was shown by Chassaing and Schaeffer that the diameter of a uniformly random quadrangulation with $n$ faces is typically of order $n^{1/4}$ \cite{cs2004rpm}.  It was then shown by Le Gall that if one rescales distances by the factor $n^{-1/4}$ \cite{lg2007top} then one obtains a tight sequence in the space of metric spaces (equipped with the Gromov-Hausdorff topology), which means that there exist subsequential limits in law.  Le Gall also showed that the Hausdorff dimension of any subsequential limit is a.s.\ equal to~$4$ \cite{lg2007top} and it was shown by Le Gall and Paulin \cite{lgp2008sphere} (see also \cite{m2008sphere}) that every subsequential limit is a.s.\ homeomorphic to the two-dimensional sphere $\s^2$.  This line of work culminated with independent works of Le Gall \cite{legalluniqueanduniversal} and Miermont \cite{miermontlimit}, which both show that one has a true limit in distribution.  The limiting random metric measure space is known as the \emph{Brownian map}.  (The term Brownian map first appeared in \cite{marckertmokkademquadrangulations}.)  We will review the continuum construction of the Brownian map in Section~\ref{subsec:tbm}.  See \cite{legall2014icm} for a more in depth survey.

Building off the works \cite{legalluniqueanduniversal,miermontlimit}, this convergence has now been extended to a number of other topologies.  In particular:
\begin{itemize}
\item Curien and Le Gall proved that the uniform quadrangulation of the whole plane (UIPQ), the local limit of a uniform quadrangulation with $n$ faces, converges to the Brownian plane \cite{cl2012brownianplane}.
\item Bettinelli and Miermont proved that quadrangulations of the disk with general boundary converge to the Brownian disk \cite{bettinelli_miermont_disks} (see also the extension \cite{gm2017simple_disk} to the case of quadrangulations of the disk with simple boundary).
\item Building on \cite{bettinelli_miermont_disks}, it was shown in \cite{gwynne-miller:uihpq,bmr2016maps} that the uniform quadrangulation of the upper half-plane (UIHPQ), the local limit of a uniform quadrangulation of the disk near a boundary typical point, converges to the Brownian half-plane.
\end{itemize}

It has long been believed that LQG should describe the large scale behavior of random planar maps (see \cite{adj1997qg}), with the case of uniformly random planar maps corresponding to $\gamma=\sqrt{8/3}$.  Precise conjectures have also been made more recently in the mathematics literature (see e.g.\ \cite{ds2011kpz,dkrv2017spheres}).

One approach is to view a quadrangulation as a surface by identifying each of the quadrilaterals with a copy of the Euclidean square $[0,1]^2$ which are identified according to boundary length using the adjacency structure of the map.  One can then conformally map the resulting surface to $\s^2$ and write the resulting metric in coordinates with respect to the Euclidean metric
\begin{equation}
\label{eqn:embedded_volume_form}
e^{\lambda_n(z)} (dx^2 + dy^2).	
\end{equation}
The goal is then to show that $\lambda_n$ converges in the limit as $n \to \infty$ to $\sqrt{8/3}$ times a form of the GFF.  A variant of this conjecture is that the volume form associated with~\eqref{eqn:embedded_volume_form} converges in the limit to the $\sqrt{8/3}$-LQG measure associated with a form of the GFF.  One can also consider other types of ``discrete'' conformal embeddings, such as circle packings, square tilings, or the Tutte (a.k.a.\ harmonic, or barycentric) embedding.  (See \cite{gms2017matedcrtmap} for a convergence result of this type in the case of the so-called \emph{mated-CRT map}.)

In this article, we will focus on the solution of a version of this conjecture carried out in \cite{qlebm,qle_continuity,qle_determined,ms2015mapmaking,quantum_spheres} in which it is shown that a $\sqrt{8/3}$-LQG surface determines a metric space structure and if one considers the correct law on $\sqrt{8/3}$-LQG surfaces then this metric space is an instance of the Brownian map.

\smallskip

\noindent{\it Acknowledgements.}  We thank Bertrand Duplantier, Ewain Gwynne, and Scott Sheffield for helpful comments on an earlier version of this article.

\section{Liouville quantum gravity surfaces}
\label{sec:surfaces}

In order to make the connections between random planar maps and LQG precise, one needs to make precise the correct law on GFF-like distributions $h$.  Since the definitions of these surfaces are quite important, we will now spend some time describing how to derive these laws (in the simply connected case), first in the setting of surfaces with boundary and then in the setting of surfaces without boundary.  These constructions were first described in \cite{she2010zipper} and carried out carefully in \cite{matingtrees}.

\subsection{Surfaces with boundary}

The starting point for deriving the correct form of the distribution $h$ is to understand the behavior of such a surface near a \emph{boundary typical point}, that is, near a point $x$ in the boundary chosen from $\nu_h$, for a general LQG surface with boundary.  To make this more concrete, we consider the domain $D = \D \cap \h$, i.e., the upper semi-disk in $\h$.  Let $h$ be a GFF on $D$ with free (resp.\ Dirichlet) boundary conditions on $[-1,1]$ (resp.\ $\partial D \setminus [-1,1]$).  Following \cite{ds2011kpz}, we then consider the law whose Radon-Nikodym derivative with respect to $h$ is given by a normalizing constant times $\nu_h([-1,1])$.  That is, if $dh$ denotes the law of $h$, then the law we are considering is given by $\nu_h([-1,1])dh$ normalized to be a probability measure.  From~\eqref{eqn:boundary_measure_eps} this, in turn, is the same as the limit as $\epsilon \to 0$ of the marginal of $h$ under the law
\begin{equation}
\label{eqn:eps_boundary_weighted}
\Theta_\epsilon(dx,dh) = \epsilon^{\gamma^2/4} e^{\gamma h_\epsilon(x)/2} dx dh	
\end{equation}
where $dx$ denotes Lebesgue measure on $[-1,1]$.  We are going to in fact consider the limit $\Theta$ as $\epsilon \to 0$ of $\Theta_\epsilon$ from~\eqref{eqn:eps_boundary_weighted}.  The reason for this is that in the limit as $\epsilon \to 0$, the conditional law of $x$ given $h$ converges to a point chosen from $\nu_h$.  By performing an integration by parts, we can write $h_\epsilon(x) = (h,\xi_\epsilon^x)_\nabla$ where $\xi_\epsilon^x(y) = -\log\max(|x-y|,\epsilon) - \wt{G}_\epsilon^x(z)$ and $\wt{G}_\epsilon^x$ is the function which is harmonic on $D$ with Dirichlet boundary conditions given by $y \mapsto -\log\max(|x-y|,\epsilon)$ on $\partial D \setminus [-1,1]$ and Neumann boundary conditions on $[-1,1]$.  In other words, $\xi_\epsilon^x$ is a truncated form of the Green's function $G$ for $\Delta$ on $D$ with Dirichlet boundary conditions on $\partial D \setminus [-1,1]$ and Neumann boundary conditions on $[-1,1]$.  Recall the following basic fact about the Gaussian distribution: if $Z \sim N(0,1)$ and we weight the law of $Z$ by a normalizing constant times $e^{\mu Z}$, then the resulting distribution is that of a $N(\mu,1)$ random variable.  By the infinite dimensional analog of this, we therefore have in the case of $h$ that weighting its law by a constant times $\exp(\gamma h_\epsilon(x)/2) = \exp(\gamma (h,\xi_\epsilon^x)_\nabla/2)$ is the same as shifting its mean by $(\gamma/2) \xi_\epsilon^x$.  That is, under $\Theta_\epsilon$, we have that the conditional law of $h$ given $x$ is that of $\wt{h} + (\gamma/2) \xi_\epsilon^x$ where $\wt{h}$ is a GFF on $D$ with free (resp.\ Dirichlet) boundary conditions on $[-1,1]$ (resp.\ $\partial D \setminus [-1,1]$).  Taking a limit as $\epsilon \to 0$, we thus see that the conditional law of $h$ given $x$ under $\Theta$ is given by that of $\wt{h} + (\gamma/2) G(x,\cdot)$ where $\wt{h}$ is a GFF on $D$ with free (resp.\ Dirichlet) boundary conditions on $[-1,1]$ (resp.\ $\partial D \setminus [-1,1]$).

This computation tells us that the local behavior of $h$ near a point chosen from $\nu_h$ is described by that of a GFF with free boundary conditions plus the singularity $-\gamma\log|\cdot|$ (as the leading order behavior of $G(x,y)$ is $-2\log|x-y|$ for $x \in (-1,1)$ and $y$ close to $x$).  We now describe how to take an infinite volume limit near $x$ in the aforementioned construction.  Roughly speaking, we will ``zoom in'' by adding a large constant $C$ to $h$ (which has the effect of replacing $\mu_h$ with $e^{\gamma C} \mu_h$), centering so that $x$ is at the origin, and then performing a rescaling so that $\D \cap \h$ is assigned one unit of mass.

It is easiest to describe this procedure if we first apply a change of coordinates from $D$ to the infinite half-strip $\strip_+ = \R_+ \times (0,\pi)$ using the unique conformal map $\varphi \colon D \to \strip_+$ which takes $-1$ to $0$, $x$ to $+\infty$, and $1$ to $\pi i$.  Then the law of $\wt{h} = h \circ \varphi^{-1} + Q \log|(\varphi^{-1})'|$ is that of $\wh{h} + (\gamma-Q) \re(\cdot)$ where $\wh{h}$ is a GFF on $\strip_+$ with free (resp.\ Dirichlet) boundary conditions on $\partial \strip_+ \setminus [0,\pi i]$ (resp.\ $[0,\pi i]$).  For each $u \geq 0$, let $\wt{A}_u$ be the average of $\wt{h}$ on the vertical line $[0,\pi i] + u$.  For such a GFF, it is possible to check that $\wt{A}_u = \wt{B}_{2u} + (\gamma-Q) u$ where $\wt{B}$ is a standard Brownian motion.  Suppose that $C > 0$ is a large constant and $\wt{h}_C = \wt{h}(\cdot+\tau_C) + C$ where $\tau_C = \inf\{u \geq 0 : \wt{A}_u + C = 0\}$.  Note that $\tau_C$ is a.s.\ finite since $(\gamma-Q) < 0$.  As $C \to \infty$, the law of $\wt{h}_C$ converges to that of a field on $\strip$ whose law can be sampled from using the following two step procedure:
\begin{itemize}
\item Take its average on vertical lines $[0,\pi i] + u$ to be given by $A_u$ where $A_u$ for $u > 0$ is given by $B_{2u} + (\gamma-Q)u$ where $B$ is a standard Brownian motion and for $u < 0$ by $\wh{B}_{-2u} + (\gamma-Q)u$ where $\wh{B}$ is an independent standard Brownian motion conditioned so that $\wh{B}_{2s} + (Q-\gamma)s \geq 0$ for all $s \geq 0$.
\item Take its projection onto the $(\cdot,\cdot)_\nabla$-orthogonal complement of the subspace of functions which are constant on vertical lines of the form $[0,\pi i] + u$ to be given by the corresponding projection of a GFF on $\strip$ with free boundary conditions on $\strip$ and which is independent of $A$.
\end{itemize}
The surface whose construction we have just described is called a \emph{$\gamma$-quantum wedge} and, when parameterized by $\h$, can be thought of as a version of the GFF on $\h$ with free boundary conditions and a $-\gamma \log|\cdot|$ singularity with the additive constant fixed in an a canonical manner.  The construction generalizes to that of an \emph{$\alpha$-quantum wedge} which is a similar type of quantum surface except with a $-\alpha \log|\cdot|$ singularity.  Quantum wedges are naturally marked by two points.  When parameterized by $\h$, these correspond to $0$ and $\infty$ and when parameterized by $\strip$ correspond to $-\infty$ and $+\infty$.  This is emphasized with the notation $(\h,h,0,\infty)$ or $(\strip,h,-\infty,+\infty)$.

A \emph{quantum disk} is the finite volume analog of a $\gamma$-quantum wedge and, when parameterized by $\strip$, the law of the associated field can be described in a manner which is analogous to that of a $\gamma$-quantum wedge.  To make this more concrete, we recall that if $B$ is a standard Brownian motion and $a \in \R$, then $e^{B_t + at}$ reparameterized to have quadratic variation $dt$ is a Bessel process of dimension $\delta = 2+2a$.  Conversely, if $Z$ is a Bessel process of dimension $\delta$, then $\log Z$ reparameterized to have quadratic variation $dt$ is a standard Brownian motion with drift $a t = (\delta-2) t/2$.  The law of the process $A$ described just above can therefore be sampled from by first sampling a Bessel process $Z$ of dimension $\delta = 8/\gamma^2$, then reparameterizing $(4/\gamma) \log Z$ to have quadratic variation $2 dt$, and then reversing and centering time so that it first hits $0$ at $u=0$.  The law of a quantum disk can be sampled from in the same manner except one replaces the Bessel process $Z$ just above with a Bessel excursion sampled from the excursion measure of a Bessel process of dimension $4-8/\gamma^2$.  A quantum disk parameterized by $\strip$ is marked by two points which correspond to $-\infty$, $+\infty$ and is emphasized with the notation $(\strip,h,-\infty,+\infty)$.  It turns out that these two points have the law of independent samples from the boundary measure when one conditions on the quantum surface structure of a quantum disk.

\subsection{Surfaces without boundary}

The derivation of the surfaces without boundary proceeds along the same lines as the derivation of the surfaces with boundary except one analyzes the behavior of an LQG surface near an area typical point rather than a boundary typical point.  The infinite volume surface is the \emph{$\gamma$-quantum cone} and the finite volume surface is the \emph{quantum sphere}.  In this case, it is natural to parameterize such a surface by the infinite cylinder $\cyl = \R \times [0,2\pi]$ (with the top and bottom identified).  The law of a $\gamma$-quantum cone parameterized by $\cyl$ can be sampled from by:
\begin{itemize}
\item Taking its average on vertical lines $[0,2\pi i] + u$ to be given by $A_u$ where $A_u$ for $u > 0$ is given by $B_u + (\gamma-Q)u$ where $B$ is a standard Brownian motion and for $u < 0$ by $\wh{B}_{-u} + (\gamma-Q)u$ where $\wh{B}$ is an independent standard Brownian motion conditioned so that $\wh{B}_s + (Q-\gamma)s \geq 0$ for all $s \geq 0$.
\item Take its projection onto the $(\cdot,\cdot)_\nabla$-orthogonal complement of the subspace of functions which are constant on vertical lines of the form $[0,2\pi i] + u$ to be given by the corresponding projection of a GFF on $\cyl$ and which is independent of $A$.
\end{itemize}
As in the case of a $\gamma$-quantum wedge, it is natural to describe $A$ in terms of a Bessel process.  In this case, one can sample from its law by first sampling a Bessel process $Z$ of dimension $\delta = 8/\gamma^2$, then reparameterizing $(2/\gamma) \log Z$ to have quadratic variation $dt$, and then reversing and centering time so that it first hits $0$ at $u=0$.  A quantum sphere can be constructed in an analogous manner except one replaces the Bessel process $Z$ just above with a Bessel excursion sampled from the excursion measure of a Bessel process of dimension $4-8/\gamma^2$.

The $\gamma$-quantum cone parameterized by $\C$ can be viewed as a whole-plane GFF plus $-\gamma \log|\cdot|$ with the additive constant fixed in a canonical way.  The $\alpha$-quantum cone is a generalization of this where the $-\gamma \log|\cdot|$ singularity is replaced with a $-\alpha \log|\cdot|$ singularity.

Quantum cones are naturally marked by two points.  When parameterized by $\C$, these correspond to $0$ and $\infty$ and when parameterized by $\cyl$ correspond to $-\infty$ and $+\infty$.  This is emphasized with the notation $(\C,h,0,\infty)$ or $(\cyl,h,-\infty,+\infty)$.  It turns out that these two points have the law of independent samples from the area measure when one conditions on the quantum surface structure of a quantum sphere.

We note that a different perspective on quantum spheres was developed in \cite{dkrv2017spheres} which follows the construction of Polyakov.  The equivalence of the construction in \cite{dkrv2017spheres} and the one described just above was established in \cite{ahs2017equiv}.  (An approach similar to \cite{ahs2017equiv} would likely yield the equivalence of the disk measures considered in \cite{lqg_disk} and the quantum disk defined earlier.)

Finally, we mention briefly that are also some works which construct LQG on non-simply connected surfaces \cite{drv2016torus,grv2016bosonic}.

\section{$\SLE$ and Liouville quantum gravity}
\label{sec:lqg_sle}

\subsection{The Schramm-Loewner evolution}
\label{subsec:sle}

The Schramm-Loewner evolution ($\SLE$) was introduced by Schramm \cite{S0} to describe the scaling limits of the interfaces in discrete lattice models in two dimensions, the motivating examples being loop-erased random walk and critical percolation.  We will discuss three variants of $\SLE$: chordal, radial, and whole-plane.

Chordal $\SLE$ is a random fractal curve which connects two boundary points in a simply connected domain.  It is most natural to define it first in $\h$ and then for other domains by conformal mapping.  Suppose that $\eta$ is a curve in $\ol{\h}$ from $0$ to $\infty$ which is non-self-crossing and non-self-tracing.  For each $t \geq 0$, let $\h_t$ be the unbounded component of $\h \setminus \eta([0,t])$ and let $g_t \colon \h_t \to \h$ be the unique conformal map with $g_t(z)-z \to 0$ as $z \to \infty$.  Then Loewner's theorem states that there exists a continuous function $W \colon [0,\infty) \to \R$ such that the maps $(g_t)$ satisfy the ODE
\begin{equation}
\label{eqn:chordal_loewner}
	\partial_t g_t(z) = \frac{2}{g_t(z) - W_t},\quad g_0(z) = z
\end{equation}
(provided $\eta$ is parameterized appropriately).  The driving function $W$ is explicitly given by $W_t = g_t(\eta(t))$.  For $\kappa > 0$, $\SLE_\kappa$ is the curve associated with the choice $W = \sqrt{\kappa} B$ where $B$ is a standard Brownian motion.  This form of the driving function arises when one makes the assumption that the law of $\eta$ is conformally invariant and satisfies the following Markov property: for each stopping time $\tau$ for $\eta$, the conditional law of $g_\tau(\eta|_{[\tau,\infty)}) - W_\tau$ is the same as the law of $\eta$.  These properties are natural to assume for the scaling limits of two-dimensional discrete lattice models.

The behavior of $\SLE_\kappa$ strongly depends on the value of $\kappa$.  When $\kappa \in (0,4]$, it describes a simple curve, when $\kappa \in (4,8)$ it is a self-intersecting curve, and when $\kappa \geq 8$ it is a space-filling curve \cite{rs2005sle}.  Special values of $\kappa$ which have been proved or conjectured to correspond to discrete lattice models include:
\begin{multicols}{2}
\begin{itemize}[leftmargin=*]
\item $\kappa=1$: Schynder woods branches \cite{lsw2017woods}
\item $\kappa=4/3$: bipolar orientation branches \cite{kmsw2015bipolar,kmsw2017bplattice}
\item $\kappa=2$: loop-erased random walk \cite{lsw2004ust}
\item $\kappa=8/3$: self-avoiding walks \cite{lsw-saw}
\item $\kappa=3$: critical Ising model \cite{smirnov-2007}
\item $\kappa=4$: level lines of the GFF \cite{ss2009contour,SchrammSheffieldGFF2}
\item $\kappa=6$: critical percolation \cite{smirnov}
\item $\kappa=16/3$: FK-Ising model \cite{smirnov-2007}
\item $\kappa=8$: uniform spanning tree \cite{lsw2004ust}
\item $\kappa=12$: bipolar orientations \cite{kmsw2015bipolar,kmsw2017bplattice}
\item $\kappa=16$: Schnyder woods \cite{lsw2017woods}
\end{itemize}
\end{multicols}

Radial $\SLE$ is a random fractal curve in a simply connected domain which connects a boundary point to an interior point.  It is defined first in $\D$ and then in other domains by conformal mapping.  The definition is analogous to the case of chordal $\SLE$, except one solves the radial Loewner ODE 
\begin{equation}
\label{eqn:radial_loewner}
	\partial_t g_t(z) = -g_t(z) \frac{g_t(z) + e^{i W_t}}{g_t(z) - e^{i W_t}}
\end{equation}
in place of~\eqref{eqn:chordal_loewner}.  As in the case of~\eqref{eqn:chordal_loewner}, the radial Loewner ODE serves to encode a non-self-crossing and non-self-tracing curve $\eta$ in terms of a continuous, real-valued function.  For each $t \geq 0$, $g_t$ is the unique conformal map from the component of $\D \setminus \eta([0,t])$ containing $0$ to $\D$ with $g_t'(0) > 0$.  Radial $\SLE_\kappa$ corresponds to the case that $W = \sqrt{\kappa} B$.  Whole-plane $\SLE$ is a random fractal curve which connects two points in the Riemann sphere.  It is defined first for the points $0$ and $\infty$ and then for other pairs of points by applying a M\"obius transformation.  It can be constructed by starting with a radial $\SLE_\kappa$ in $\C \setminus (\epsilon \D)$ from $\epsilon$ to $\infty$ and then taking a limit as $\epsilon \to 0$.

\subsection{Exploring an LQG surface with an $\SLE$}
\label{subsec:exploring_lqg}

There are two natural operations that one can perform in the context of planar maps (see the physics references in \cite{ds2011kpz}).  Namely:
\begin{itemize}
\item One can ``glue'' together two planar maps with boundary by identifying their edges along a marked boundary segment to produce a planar map decorated with a distinguished interface.  If the two maps are chosen independently and uniformly at random, then this interface will in fact be a self-avoiding walk (SAW).  (See Section~\ref{subsec:saw} for more details.)
\item One can also decorate a planar map with an instance of a statistical physics model (e.g.\ a critical percolation configuration or a uniform spanning tree) and then explore the interfaces of the statistical physics model.  (See Section~\ref{subsec:perc} for more details in the case of percolation.)
\end{itemize}

If one takes as an ansatz that LQG describes the large scale behavior of random planar maps, then it is natural to guess that one should be able to perform the same operations in the continuum on LQG.  In order to make this mathematically precise, one needs to describe the precise form of the laws of:
\begin{itemize}
\item The field $h$ which describes the underlying LQG surface and
\item The law of the interfaces.
\end{itemize}
In view of the discussion in Section~\ref{sec:surfaces}, it is natural to expect that quantum wedges, disks, cones, and spheres will play the role of the former.  In view of the conformal invariance ansatz for critical models in two-dimensional statistical mechanics, it is natural to expect that $\SLE$-type curves should play the role of the latter.

Recall that LQG comes with the parameter~$\gamma$ and $\SLE$ comes with the parameter~$\kappa$.  As we will describe below in more detail, it is important that these parameters are correctly tuned.  Namely, it will always be the case that
\begin{equation}
\label{eqn:gamma_kappa}
\gamma = \min\left( \sqrt{\kappa}, \frac{4}{\sqrt{\kappa}} \right).	
\end{equation}
We emphasize that~\eqref{eqn:gamma_kappa} states that for each $\gamma \in (0,2)$, there are precisely two compatible values of $\kappa$: $\kappa = \gamma^2 \in (0,4)$ and $\kappa = 16/\gamma^2 > 4$.

We will now describe some work of Sheffield \cite{she2010zipper}, which is the first mathematical result relating $\SLE$ to LQG and should be interpreted as the continuous analog of the gluing operation for planar maps with boundary.  It is motivated by earlier work of Duplantier from the physics literature \cite{dup1998rw_qg,duplantier1999percolation,duplantier1999saw,dup2000fractals}.

\begin{theorem}
\label{thm:zipper}
Fix $\kappa \in (0,4)$ and let $\gamma = \sqrt{\kappa}$.  Suppose that $\cW = (\h,h,0,\infty)$ is a  $(\gamma-2/\gamma)$-quantum wedge and that $\eta$ is an independent $\SLE_\kappa$ process in $\h$ from $0$ to $\infty$.  Let $D_1$ (resp.\ $D_2$) be the component of $\h \setminus \eta$ which is to the left (resp.\ right) of $\eta$.  Then the quantum surfaces $\cD_1 = (D_1,h,0,\infty)$ and $\cD_2 = (D_2,h,0,\infty)$ are independent $\gamma$-quantum wedges.  Moreover, $\cW$ and $\eta$ are a.s.\ determined by $\cD_1,\cD_2$.
\end{theorem}

We first emphasize that the independence of $\cD_1, \cD_2$ in Theorem~\ref{thm:zipper} is in terms of quantum surfaces, which are themselves defined modulo conformal transformation.  The $\gamma$-quantum wedges $\cD_1,\cD_2$ which are parameterized by the regions which are to the left and right of $\eta$ each have their own boundary length measure.  Therefore for any point~$z$ along~$\eta$, one can measure the boundary length distance from~$z$ to~$0$ along the left or the right side of~$\eta$ (i.e., using $\cD_1$ or $\cD_2$).  One of the other main results of \cite{she2010zipper} is that these two quantities agree.  An $\SLE_\kappa$ curve with $\kappa \in (0,4)$ therefore has a well-defined notion of quantum length.  Moreover, this also allows one to think of \cite{she2010zipper} as a statement about welding quantum surfaces together.

The cutting/welding operation first established in \cite{she2010zipper} was substantially generalized in \cite{matingtrees}, in which many other $\SLE$ explorations of LQG surfaces were studied.  Let us mention one result in the context of an $\SLE_\kappa$ process for $\kappa \in (4,8)$.

\begin{theorem}
\label{thm:sle6}
Fix $\kappa \in (4,8)$ and let $\gamma = 4/\sqrt{\kappa}$.  Suppose that $\cW = (\h,h,0,\infty)$ is a $(4/\gamma-\gamma/2)$-quantum wedge and that $\eta$ is an independent $\SLE_\kappa$ process.  Then the quantum surfaces parameterized by the components of $\h \setminus \eta$ are conditionally independent quantum disks given their boundary lengths.  The boundary lengths of these disks which are on the left (resp.\ right) side of $\eta$ are in correspondence with the jumps of a $\kappa/4$-stable L\'evy process $L$ (resp.\ $R$) with only downward jumps.  Moreover, $L$ and $R$ are independent.
\end{theorem}

The time parameterization of the L\'evy processes $L$ and $R$ gives rise to an intrinsic notion of time for $\eta$ which in \cite{matingtrees} is referred to as the \emph{quantum natural time}.

There are also similar results for explorations of the finite volume surfaces by $\SLE$ processes.  We will focus on one such result here (proved in \cite{quantum_spheres}) which is relevant for the construction of the metric on $\sqrt{8/3}$-LQG.

\begin{theorem}
\label{thm:sphere}
Suppose that $\gamma=\sqrt{8/3}$ and that $\cS = (\cyl,h,-\infty,+\infty)$ is a quantum sphere.  Let $\eta$ be a whole-plane $\SLE_6$ process in $\cyl$ from $-\infty$ to $+\infty$ which is sampled independently of $h$ and then reparameterized according to quantum natural time.  Let $X_t$ be the quantum boundary length of the component $C_t$ of $\cyl \setminus \eta([0,t])$ containing $+\infty$.  Then $X$ evolves as the time-reversal of a $3/2$-stable L\'evy excursion with only upward jumps.  For a given time $t$, the conditional law of the surface parameterized by $C_t$ given $X_t$ is that of a quantum disk with boundary length $X_t$ weighted by its quantum area and the conditional law of the point $\eta(t)$ is given by the quantum boundary length measure on $\partial C_t$.  Moreover, the surfaces parameterized by the other components of $\cyl \setminus \eta([0,t])$ are quantum disks which are conditionally independent given $X|_{[0,t]}$, each correspond to a downward jump of $X|_{[0,t]}$, and have quantum boundary length given by the corresponding jump.
\end{theorem}

\section{Construction of the metric}

In this section we will describe the construction of the metric on $\sqrt{8/3}$-LQG from \cite{qlebm,qle_continuity,qle_determined}.  The construction is strongly motivated by discrete considerations, which we will review in Section~\ref{subsec:eden}, before reviewing the continuum construction in Section~\ref{subsec:qle}.

\subsection{The Eden growth model}
\label{subsec:eden}

Suppose that $G = (V,E)$ is a connected graph.  In the Eden growth model \cite{eden1961two} (or first-passage percolation \cite{hammersley1965first}), one associates with each $e \in E$ an independent $\exp(1)$ weight $Z_e$.  The aim is then to understand the \emph{random} metric $d_{\mathrm{FPP}}$ on $G$ which assigns length $Z_e$ to each $e \in E$.  Suppose that $x \in V$.  By the memoryless property of the exponential distribution, there is a simple Markovian way of growing the $d_{\mathrm{FPP}}$-metric ball centered at $x$.  Namely, one inductively defines an increasing sequence of clusters $C_n$ as follows.
\begin{itemize}
\item Set $C_0 = \{x\}$.
\item Given that $C_n$ is defined, choose an edge $e = \{y,z\}$ uniformly at random among those with $y \in C_n$ and $z \notin C_n$, and then take $C_{n+1} = C_n \cup \{z\}$.
\end{itemize}

\begin{figure}[ht!]
\begin{center}
\includegraphics[width=0.45\textwidth]{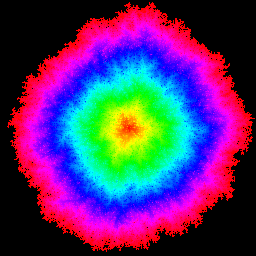} \hspace{0.03\textwidth}
\includegraphics[width=0.45\textwidth]{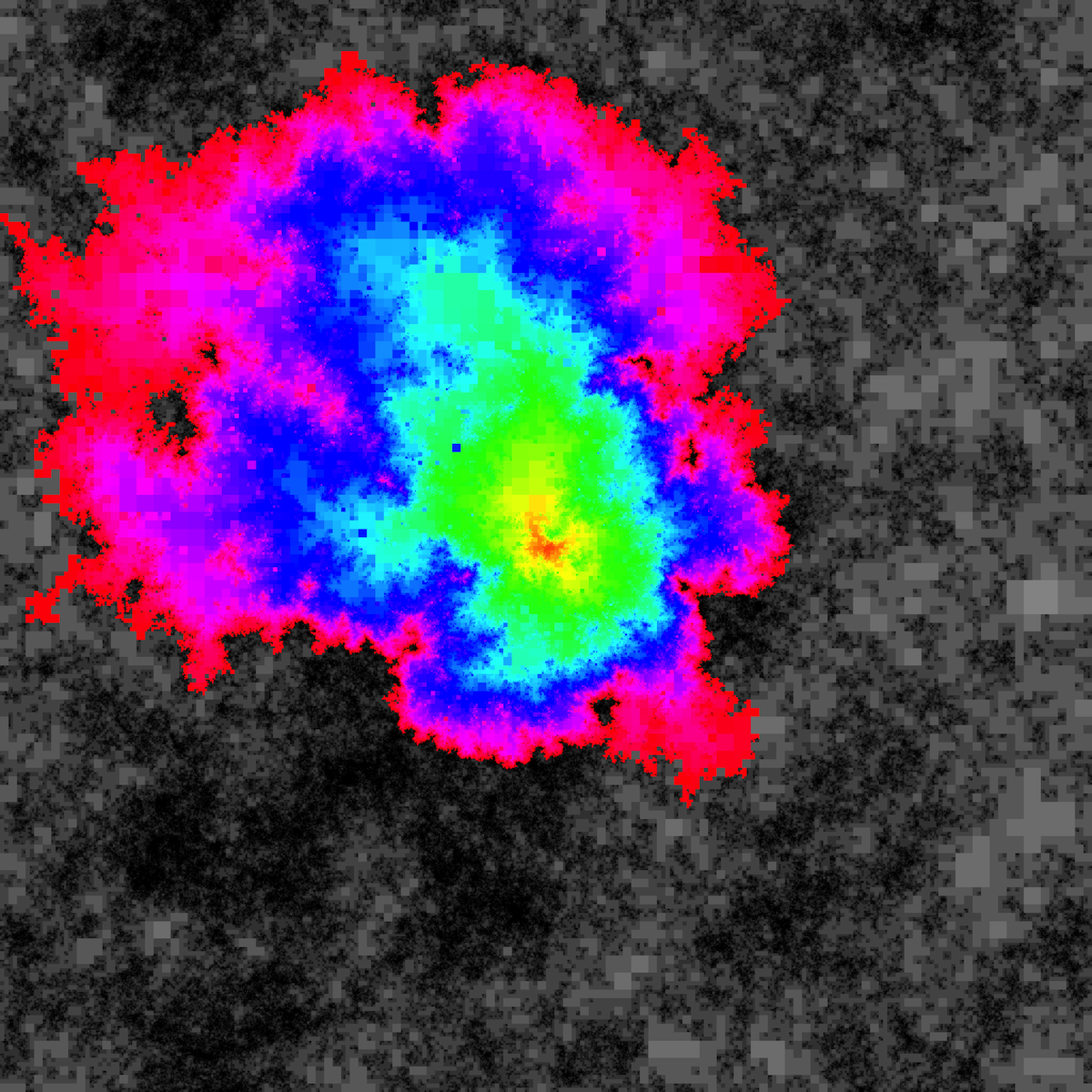}
\end{center}
\caption{\label{fig:eden}{\bf Left:} Eden model on $\Z^2$.  {\bf Right:} Eden model on a graph approximation of $\sqrt{8/3}$-LQG.  This serves as a discretization of $\QLE(8/3,0)$.}
\end{figure}

It is an interesting question to analyze the large-scale behavior of the clusters $C_n$ on a given graph $G$.  One of the most famous examples is the case $G = \Z^2$, i.e., the two-dimensional integer lattice (see the left side of Figure~\ref{fig:eden}).  It was shown by Cox and Durrett \cite{cd1981fpp} that the macroscopic shape of $C_n$ is convex but computer simulations suggest that it is not a Euclidean ball.  The reason for this is that $\Z^2$ is not sufficiently isotropic.  In \cite{vw1990fpp,vw1992fpp}, Wiermann and Vahidi-Asl considered the Eden model on the Delaunay triangulation associated with the Voronoi tesselation of a Poisson point process with Lebesgue intensity on $\R^2$ and showed that at large scales the Eden growth model is well-approximated by a Euclidean ball.  The reason for the difference between the setting considered in \cite{vw1990fpp,vw1992fpp} and $\Z^2$ is that such a Poisson point process does not have preferential directions due to the underlying randomness and the rotational invariance of Lebesgue measure.

Also being random, it is natural to expect that the Eden growth model on a random planar map should at large scales be approximated by a metric ball.  This has now been proved by Curien and Le Gall \cite{clg2015distance} in the case of the planar dual of a random triangulation.  The Eden growth model on the planar dual of a random triangulation is natural to study because it can be described in terms of a so-called \emph{peeling process}, in which one picks an edge uniformly on the boundary of the cluster so far and then reveals the opposing triangle.  In particular, this exploration procedure respects the Markovian structure of a random triangulation in the sense that one can describe the law of the regions cut out by the growth process (uniform triangulations of the disk) as well as the evolution of the boundary length of the cluster.  We will eventually want to make sense of the Eden model on a quantum sphere which satisfies the same properties.  In order to motivate the construction, we will describe a variant of the Eden growth model on the planar dual of a random triangulation which involves two operations we know how to perform in the continuum (in view of the connection between $\SLE$ and LQG) and respects the Markovian structure of the quantum sphere in the same way as described above.  Namely, we fix $k \in \N$ and then define an increasing sequence of clusters $C_n$ in the dual of the map as follows.
\begin{itemize}
\item Set $C_0 = \{F\}$ where $F$ is the root face.
\item Given that $C_n$ is defined, pick two edges $e_1,e_2$ on the outer boundary of $C_n$ uniformly at random and color the vertices on the clockwise (resp.\ counterclockwise) arc of the outer boundary of $C_n$ from $e_1$ to $e_2$ blue (resp.\ yellow).
\item Color the vertices in the remainder of the map blue (resp.\ yellow) independently with probability $1/2$.  Take $C_{n+1}$ to be the union of $C_n$ and $k$ edges which lie along the blue/yellow percolation interface starting from one of the marked edges described above.	
\end{itemize}
This growth process can also be described in terms of a peeling process, hence it respects the Markovian structure of a random triangulation in the same way as the Eden growth model.  It is also natural to expect there to be universality.  Namely, the macroscopic behavior of the cluster should not depend on the specific geometry of the ``chunks'' that we are adding at each stage.  That is, it should be the case that for each fixed $k \in \N$, the cluster $C_n$ at large scales should be well-approximated by a metric ball.

\subsection{$\QLE(8/3,0)$: The Eden model on the quantum sphere}
\label{subsec:qle}

Suppose that $(\cS,x,y)$ is a quantum sphere marked by two independent samples $x,y$ from the quantum measure.  We will now describe a variant of the above construction on $\cS$ starting from $x$ and targeted at $y$, with $\SLE_6$ as the continuum analog of the percolation interface.  Fix $\delta > 0$ and let $\eta_0$ be a whole-plane $\SLE_6$ on $\cS$ from $x$ to $y$ with the quantum natural time parameterization as in Theorem~\ref{thm:sle6}.  Then Theorem~\ref{thm:sle6} implies that:
\begin{itemize}
\item The quantum surfaces parameterized by the components of $\cS \setminus \eta_0([0,\delta])$ which do not contain $y$ are conditionally independent quantum disks.
\item The component $C_0$ of $\cS \setminus \eta_0([0,\delta])$ containing $y$ has the law of a quantum disk weighted by its quantum area.
\item The conditional law of $\eta_0(\delta)$ is given by the quantum boundary measure on $\partial C_0$.
\end{itemize}
We now pick $z_1$ from the quantum boundary measure on $\partial C_0$ and then let $\eta_1$ be a radial $\SLE_6$ in $C_0$ from $z_1$ to $y$, parameterized by quantum natural time.  Then it also holds that the quantum surfaces parameterized by the components of $C_0 \setminus \eta_1([0,\delta])$ which do not contain $y$ are conditionally independent quantum disks, the component $C_1$ of $C_0 \setminus \eta_1([0,\delta])$ containing $y$ has the law of a quantum disk weighted by its quantum area, and $\eta_1(\delta)$ is uniformly distributed according to the quantum boundary measure on $\partial C_1$.

The $\delta$-approximation $\Gamma^{\delta,x\to y}$ to $\QLE(8/3,0)$ (\emph{quantum Loewner evolution} with parameters $\gamma^2=8/3$ and $\eta=0$; see Section~\ref{subsec:qle_general} for more on the family of processes $\QLE(\gamma^2,\eta)$) from $x$ to $y$ is defined by iterating the above procedure until it eventually reaches $y$.  One can view $\Gamma^{\delta,x\to y}$ as arising by starting with a whole-plane $\SLE_6$ process $\eta$ on $\cS$ from $x$ to $y$ parameterized by quantum natural time and then resampling the location of the tip of $\eta$ at each time of the form $k \delta$ where $k \in \N$.  Due to the way the process is constructed, we emphasize that the following hold:
\begin{itemize}
\item The surfaces parameterized by the components of $\cS \setminus \Gamma_t^{\delta,x \to y}$ which do not contain $y$ are conditionally independent quantum disks.
\item The surface parameterized by the component $C_t$ of $\cS \setminus \Gamma_t^{\delta,x \to y}$ which contains $y$ on its boundary is a quantum disk weighted by its area.
\item The evolution of the quantum boundary length $X_t$ of $\partial C_t$ is the same as in the case of $\eta$, i.e., it is given by the time-reversal of a $3/2$-stable L\'evy excursion.
\end{itemize}
That is, $\Gamma^{\delta,x \to y}$ respects the Markovian structure of a quantum sphere.  The growth process $\QLE(8/3,0)$ is constructed by taking a limit as $\delta \to 0$ of the $\delta$-approximation defined above.  All of the above properties are preserved by the limit.

Parameterizing an $\SLE_6$ by quantum natural time on a quantum sphere is the continuum analog of parameterizing a percolation exploration on a random planar map according to the number of edges that the percolation has visited.  This is not the correct notion of time if one wants to define a metric space structure since one should be adding particles to the growth process at a rate which is proportional to its boundary length.  One is thus led to make the following change of time: set
\[ D(t) = \int_0^t \frac{1}{X_u} du\]
and then let $s(r) = \inf\{t \geq 0 : D(t) > r\}$.  Due to the above interpretation, the time-parameterization $s(r)$ is called the \emph{quantum distance time}.  Let $\Gamma^{x \to y}$ be a $\QLE(8/3,0)$ parameterized by $s(r)$ time.  Then it will ultimately be the case that $\Gamma_r^{x \to y}$ defines a metric ball of radius $r$, and we will sketch the proof of this fact in what follows.

\newcommand{\qdist}{d_{\cQ}}

Let $(x_n)$ be a sequence of i.i.d.\ points chosen from the quantum measure on $\cS$.  For each $i \neq j$, we let $\Gamma^{x_i \to x_j}$ be a conditionally independent (given $\cS$) $\QLE(8/3,0)$ from $x_i$ to $x_j$ with the quantum distance parameterization.  We then set $\qdist(x_i,x_j)$ be to be the amount of time it takes for $\Gamma^{x_i \to x_j}$ to reach $x_j$.  We want to show that $\qdist(x_i,x_j)$ defines a metric on the set $(x_i)$ which is a.s.\ determined by $\cS$.  It is obvious from the construction that $\qdist(x_i,x_j) > 0$ for $i \neq j$, so to establish the metric property it suffices to prove that $\qdist$ is symmetric and satisfies the triangle inequality.  This is proved by making use of a strategy developed by Sheffield, Watson, and Wu in the context of $\CLE_4$.

Symmetry, the triangle inequality, and the fact that~$\qdist$ is a.s.\ determined by~$\cS$ all follow from the following stronger statement (taking without loss of generality $x=x_1$ and $y=x_2$).  Let $\Theta$ denote the law of $(\cS,x,y,\Gamma^{x \to y}, \Gamma^{y \to x},U)$ where $U$ is uniform in $[0,1]$ independently of everything else.  Let $\Theta^{x \to y}$ (resp.\ $\Theta^{y \to x}$) be the law whose Radon-Nikodym derivative with respect to $\Theta$ is given by $\qdist(x,y)$ (resp.\ $\qdist(y,x)$).  That is,
\[ \frac{d\Theta^{x \to y}}{d\Theta} = \qdist(x,y) \quad\text{and}\quad \frac{d\Theta^{y \to x}}{d\Theta} = \qdist(y,x).\]
We want to show that $\Theta^{x \to y} = \Theta^{y \to x}$ because then the uniqueness of Radon-Nikodym derivatives implies that $\qdist(x,y) = \qdist(y,x)$.  Since $\Gamma^{x \to y}$ and $\Gamma^{y \to x}$ were taken to be conditionally independent given $\cS$, this also implies that the common value of $\qdist(x,y)$ and $\qdist(y,x)$ is a.s.\ determined by $\cS$.

The main step in proving this is the following, which is a restatement of \cite[Lemma~1.2]{qlebm}.
\begin{lemma}
\label{lem:symmetry}
Let $\tau = U \qdist(x,y)$ so that $\tau$ is uniform in $[0,\qdist(x,y)]$ and let $\ol{\tau} = \inf\{t \geq 0 : \Gamma_\tau^{x \to y} \cap \Gamma_t^{y \to x} \neq \emptyset\}$.  We similarly let $\ol{\sigma} = U \qdist(y,x)$ and $\sigma = \inf\{t \geq 0 : \Gamma_t^{x \to y} \cap \Gamma_{\ol{\sigma}}^{y \to x} \neq \emptyset\}$.  Then the $\Theta^{x \to y}$ law of $(\cS,x,y,\Gamma^{x \to y}|_{[0,\tau]}, \Gamma^{y \to x}|_{[0,\ol{\tau}]})$ is the same as the $\Theta^{y \to x}$ law of $(\cS,x,y,\Gamma^{x \to y}|_{[0,\sigma]},\Gamma^{y \to x}|_{[0,\ol{\sigma}]})$.
\end{lemma}
 Upon proving Lemma~\ref{lem:symmetry}, the proof is completed by showing that the $\Theta^{x \to y}$ conditional law of $\Gamma^{x \to y}, \Gamma^{y \to x}$ given $(\cS,x,y,\Gamma^{x \to y}|_{[0,\tau]}, \Gamma^{y \to x}|_{[0,\ol{\tau}]})$ is the same as the $\Theta^{y \to x}$ conditional law of $\Gamma^{x \to y}, \Gamma^{y \to x}$ given $(\cS,x,y,\Gamma^{x \to y}|_{[0,\sigma]},\Gamma^{y \to x}|_{[0,\ol{\sigma}]})$.

Lemma~\ref{lem:symmetry} turns out to be a consequence of a corresponding symmetry statement for whole-plane $\SLE_6$, which in a certain sense reduces to the time-reversal symmetry of whole-plane $\SLE_6$ \cite{ms2013imag4}, and then reshuffling the $\SLE_6$ as described above to obtain a $\QLE(8/3,0)$.

The rest of the program carried out in \cite{qle_continuity,qle_determined} consists of showing that:
\begin{itemize}
\item The metric defined above extends uniquely in a continuous manner to the entire quantum sphere yielding a metric space which is homeomorphic to~$\s^2$.  Moreover, the resulting metric space is geodesic and isometric to the Brownian map using the characterization from \cite{ms2015mapmaking}.  We will describe this in further detail in the next section.
\item The quantum sphere instance is a.s.\ determined by the metric measure space structure.  This implies that the Brownian map possesses a canonical embedding into $\s^2$ which takes it to a form of $\sqrt{8/3}$-LQG.  This is based on an argument which is similar to that given in \cite[Section~10]{matingtrees}.
\end{itemize}

\subsection{$\QLE(\gamma^2,\eta)$}
\label{subsec:qle_general}

\begin{figure}[ht!]
\begin{center}
\includegraphics[width=0.31\textwidth]{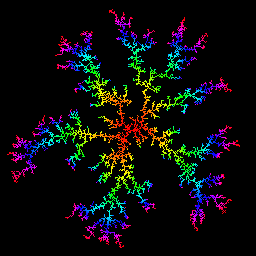} \hspace{0.005\textwidth}
\includegraphics[width=0.31\textwidth]{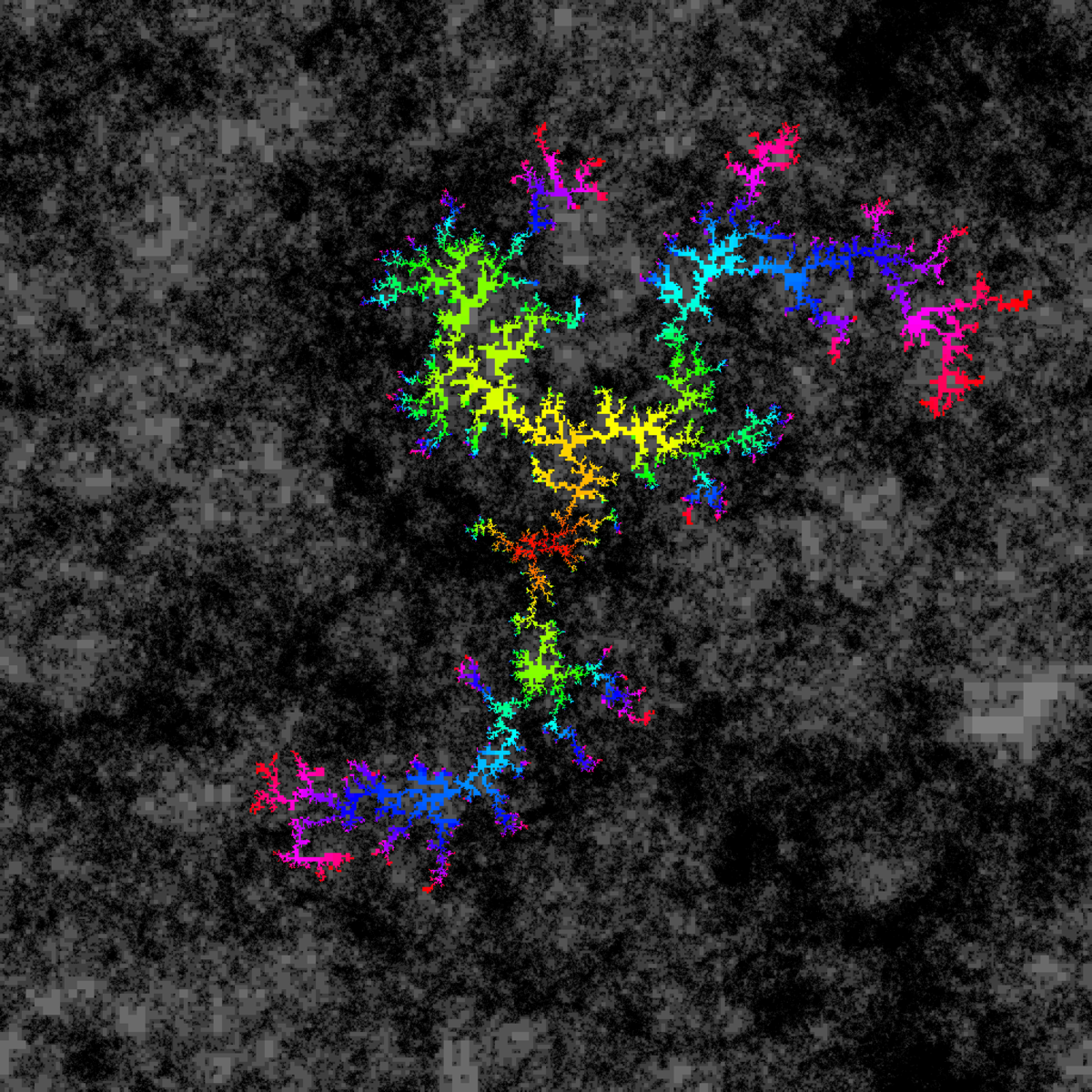}
\hspace{0.005\textwidth}
\includegraphics[width=0.33\textwidth]{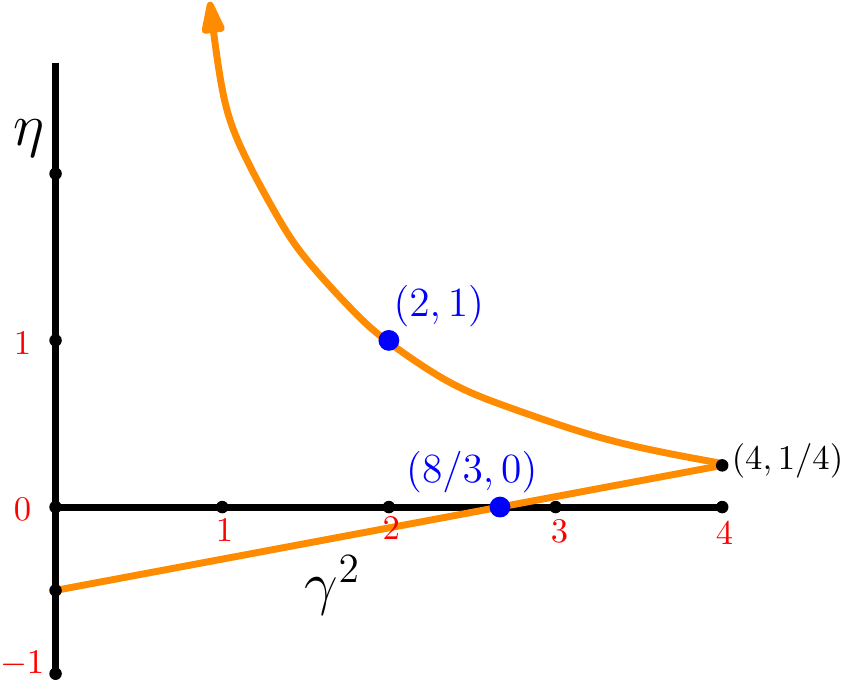}
\end{center}
\caption{\label{fig:dla_qle} {\bf Left:} DLA on $\Z^2$.  {\bf Middle:} DLA on a graph approximation of $\sqrt{2}$-LQG.  This serves as a discretization of $\QLE(2,1)$.  {\bf Right:} Plot of $(\gamma^2,\eta)$ values for which $\QLE(\gamma^2,\eta)$ is constructed in \cite{ms2013qle}.}
\end{figure}
The process $\QLE(8/3,0)$ described in Section~\ref{subsec:qle} is a part of a general family of growth processes which are the conjectural scaling limits of a family of growth models on random surfaces which we now describe.

The dielectric breakdown model (DBM) \cite{niemeyer1984fractal} with parameter $\eta$ is a family of models which interpolate between the Eden model and diffusion limited aggregation (DLA).  If $\mu_{{\rm HARM}}$ (resp.\ $\mu_{{\rm LEN}}$) denotes the natural harmonic (resp.\ length) measure on the underlying surface, then microscopic particles are added according to the measure
\[ \left( \frac{d\mu_{{\rm HARM}}}{d \mu_{{\rm LEN}}} \right)^\eta d\mu_{{\rm LEN}}.\]
In particular, $\eta = 0$ corresponds to the Eden model and $\eta=1$ corresponds to DLA.

One can apply the tip-rerandomization procedure to $\SLE_\kappa$ or $\SLE_{\kappa'}$ coupled with $\gamma$-LQG for other values of $\kappa,\kappa'$ and $\gamma$ provided $\kappa = \gamma^2$ and $\kappa'=16/\gamma^2$.  The resulting process, which is called $\QLE(\gamma^2,\eta)$ and is defined and analyzed in \cite{ms2013qle}, is the conjectural scaling limit of $\eta$-DBM on a $\gamma$-LQG surface where
\[ \eta = \frac{3\gamma^2}{16} - \frac{1}{2} \quad\text{or} \quad \eta = \frac{3}{\gamma^2} - \frac{1}{2}.\]
(The parameter $\gamma$ determines the type of LQG surface on which the process grows and the value of $\eta$ determines the manner in which it grows.)  Special parameter values include $\QLE(8/3,0)$ and $\QLE(2,1)$, the latter being the conjectural scaling limit of DLA on a tree-weighted random planar map.  See Figure~\ref{fig:dla_qle}.  It remains an open question to construct $\QLE(\gamma^2,\eta)$ for the full range of parameter values.

\section{Equivalence with the Brownian map}

In the previous section, we have described how to build a metric space structure on top of $\sqrt{8/3}$-LQG.  We will describe here how it is checked that this metric space structure is equivalent with the Brownian map.

\subsection{Brownian map definition}
\label{subsec:tbm}

The Brownian map is constructed from a continuum version of the Cori-Vauquelin-Schaeffer bijection \cite{cv1981maps,schaeffer} using the \emph{Brownian snake}.  Suppose that $Y \colon [0,T] \to \R_+$ is picked from the excursion measure for Brownian motion.  Given $Y$, we let $X$ be a centered Gaussian process with $X_0 = 0$ and
\[ \cov(X_s,X_t) = \inf\{ Y_r : r \in [s,t]\}.\]
For $s < t$ and $[t,s] = [0,T] \setminus (s,t)$, we set
\[ d^\circ(s,t) = X_s + X_t - 2\max\left( \min_{r \in [s,t]} X_r, \min_{r \in [t,s]} X_r \right).\]
Let $\cT$ be the instance of the continuum random tree (CRT) \cite{ald1991crt1} encoded by $Y$ and let $\rho \colon [0,T] \to \cT$ be the corresponding projection map.  We then set
\[ d_\cT(a,b) = \min\{ d^\circ(s,t) : \rho(s) = a,\quad \rho(t) = b\}.\]
Finally, for $a,b \in \cT$, we set
\[ d(a,b) = \inf\left\{ \sum_{j=1}^k d_\cT^\circ(a_{j-1},a_j) \right\}\]
where the infimum is over all $k \in \N$ and $a_0 = a,a_1,\ldots,a_k = b$ in $\cT$.  Quotienting by the equivalence relation $a \cong b$ if and only if $d(a,b) = 0$ yields a metric space $(S,d)$.   It is naturally equipped with a measure $\nu$ by taking the projection of Lebesgue measure on $[0,T]$.  Finally, $(S,d,\nu)$ is naturally marked by the points $x$ and $y$ which are respectively given by the projections of $t=0$ and the value of $t$ at which $X$ attains its infimum.  The space $(S,d,\nu,x,y)$ is the (doubly marked) Brownian map and we denote its law by $\mustwo$.  It is an infinite measure (since the Brownian excursion measure is an infinite measure).  As mentioned earlier, it follows from \cite{lgp2008sphere} (see also \cite{m2008sphere}) that $(S,d)$ a.s.\ has the topology of $\s^2$.  Also, the law of $(S,d,\nu,x,y)$ is invariant under the operation of resampling $x,y$ independently from $\nu$.  The standard unit area Brownian map arises by conditioning $\mustwo$ to have total area equal to $1$.  Equivalently, one can take the construction above and condition on the Brownian excursion $Y$ to have length equal to $1$.

Variants of the Brownian map with other topologies are defined in a similar manner.  For example, the Brownian disk, half-plane, and plane are defined this way in \cite{bettinelli_miermont_disks,alg,lg_disk_snake,bmr2016maps,gwynne-miller:uihpq,cl2012brownianplane}.

\subsection{The $\alpha$-stable L\'evy net}

Suppose that we have a doubly-marked metric space $(S,d,x,y)$ which has the topology of $\s^2$.  For each $r \in [0,d(x,y)]$, we define the \emph{filled metric ball} $\fb{x}{r}$ to be the closure of the complement of the $y$-containing component of $S \setminus B(x,r)$.  The \emph{metric net} of $(S,d,x,y)$ is the closure of the union of $\partial \fb{x}{r}$ over $r \in [0,d(x,y)]$.  It turns out that it is possible to give an explicit description of the law of the metric net of the Brownian map and, as we will explain just below, this is one of the main ingredients which characterizes its law.

\begin{figure}[ht!]
\begin{center}
\includegraphics[scale=.7]{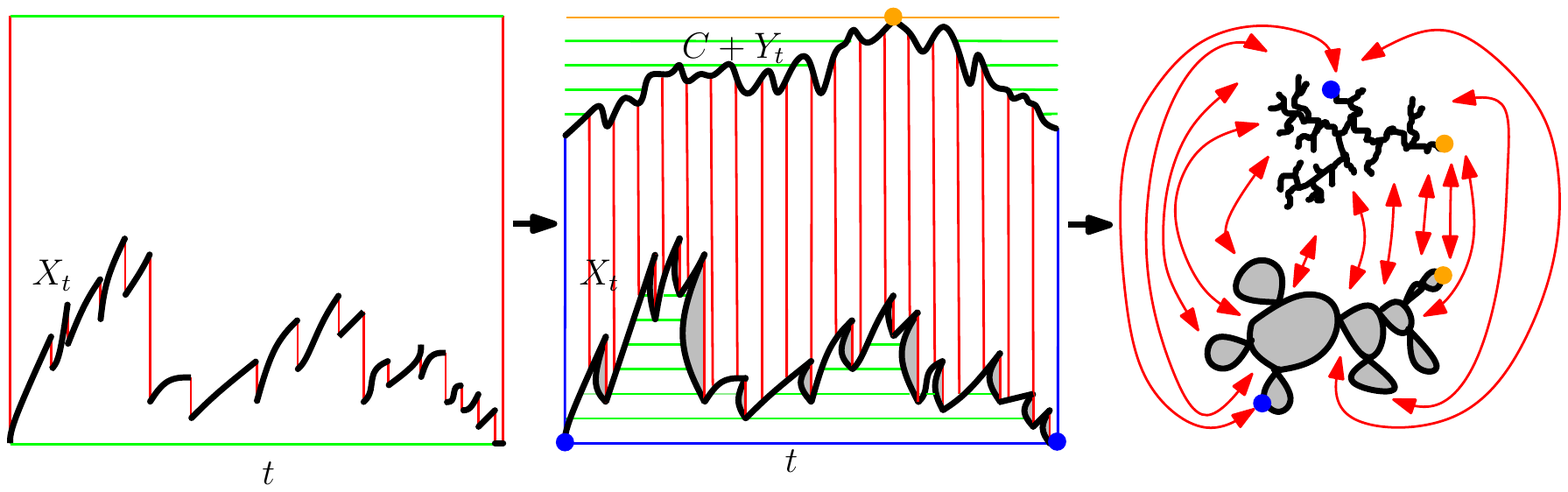}
\caption{\label{fig:levynet} Illustration of the construction of the $\alpha$-stable L\'evy net.}
\end{center}
\end{figure}

The \emph{$\alpha$-stable L\'evy net} is defined as follows.  Fix $\alpha \in (1,2)$ and suppose that $X$ is the time-reversal of an $\alpha$-stable L\'evy excursion with only upward jumps (so that $X$ has only downward jumps).  We define the height process $Y$ associated with $X$ to be given at time $t$ by the amount of local time that $X|_{[t,T]}$ spends at its running infimum. Both $X$ and $Y$ encode trees.  Namely, associated with $X$ is a looptree which is obtained by considering the graph of $X$ and then replacing each of the downward jumps by a topological disk that does not otherwise cross the graph of $X$ and with the disks pairwise disjoint.  Points on the graph of~$X$ or the disk boundaries are considered to be equivalent if they can be connected by a horizontal chord which lies entirely below the graph of $X$.  The tree associated with~$Y$ is defined by declaring two points on the graph of~$Y$ to be equivalent if they can be connected by a horizontal chord which lies entirely above the graph of~$Y$.  We can then glue these two trees together as illustrated in Figure~\ref{fig:levynet} to obtain the $\alpha$-stable L\'evy net.  The tree associated with~$X$ (resp.\ $Y$) is the called the dual (resp.\ geodesic) tree of the L\'evy net instance.

Due to the construction, points on the geodesic tree which are identified with each other all have the same distance to the root in the geodesic tree.  Therefore every point in the L\'evy net has a well-defined distance to the root, hence one can talk about metric balls which are centered at the root.  It is not difficult to see that one can in fact associate with each such metric ball a boundary length and that this boundary length evolves as a so-called continuous state branching process (CSBP) as one reduces the radius of the ball.  In fact, the boundary lengths between any finite collection of geodesics evolve as collection of independent CSBPs as the ball radius is reduced and this property essentially characterizes the $\alpha$-stable L\'evy net (as it determines how the geodesics are glued together).

It is proved in \cite{ms2015mapmaking} that the law of the metric net of a sample from $\mustwo$ is given by that of a $3/2$-stable L\'evy net.  The following is a restatement of \cite[Theorem~4.6]{ms2015mapmaking}.

\begin{theorem}
\label{thm:characterization}
The doubly marked Brownian map measure $\mustwo$ is the unique infinite measure on doubly-marked metric measure spaces $(S,d,\nu,x,y)$ with the topology of $\s^2$ which satisfy the following properties:
\begin{enumerate}
\item The law of $(S,d,\nu,x,y)$ is invariant under resampling $x$ and $y$ independently from $\nu$.
\item The law of the metric net from $x$ to $y$ agrees with that of the $\alpha$-stable L\'evy net for some value of $\alpha \in (1,2)$.
\item For each fixed $r > 0$, the metric measure spaces $\fb{x}{r}$ and $S \setminus \fb{x}{r}$ (equipped with the interior internal metric and the restriction of $\nu$) are conditionally independent given the boundary length of $\partial \fb{x}{r}$.
\end{enumerate}
\end{theorem}

The proof of Theorem~\ref{thm:characterization} given in \cite{ms2015mapmaking} shows that the assumptions necessarily imply that $\alpha=3/2$, which is why we do not need to make this assumption explicitly in the statement of Theorem~\ref{thm:characterization}.

\subsection{$\QLE(8/3,0)$ metric satisfies the axioms which characterize the Brownian map}

To check that the $\QLE(8/3,0)$ metric defined on the $\sqrt{8/3}$-quantum sphere defines an instance of the Brownian map, it suffices to check that the axioms of Theorem~\ref{thm:characterization} are satisfied.  The construction of the $\QLE(8/3,0)$ metric in fact implies that the first and third axioms are satisfied, so the main challenge is to check the second axiom.  This is one of the aims of \cite{qle_continuity} and is closely related to the form of the evolution of the boundary length when one performs an $\SLE_6$ exploration on a quantum sphere as described in Theorem~\ref{thm:sphere}.

Upon proving the equivalence of the $\QLE(8/3,0)$ metric on a quantum sphere with the Brownian map, it readily follows that several other types of quantum surfaces with $\gamma=\sqrt{8/3}$ are equivalent to certain types of Brownian surfaces.  Namely, the Brownian disk, half-plane, and plane are respectively equivalent to the quantum disk, $\sqrt{8/3}$-quantum wedge, and $\sqrt{8/3}$-quantum cone \cite{qle_continuity,gwynne-miller:uihpq,lg_disk_snake}.

\section{Scaling limits}

The equivalence of LQG and Brownian surfaces allows one to define $\SLE$ on a Brownian surface in a canonical way as the embedding of a Brownian surface is a.s.\ determined by the metric measure space structure \cite{qle_determined}.  This makes it possible to prove that certain statistical physics models on uniformly random planar maps converge to $\SLE$.  The natural topology of convergence is the so-called Gromov-Hausdorff-Prokhorov-uniform topology (GHPU) developed in \cite{gwynne-miller:uihpq}, which is an extension of the Gromov-Hausdorff topology to curve-decorated metric measure spaces.  (We note that a number of other scaling limit results for random planar maps decorated with a statistical mechanics model toward $\SLE$ on LQG have been proved in e.g.\ \cite{sheffield2011qg_inventory,kmsw2015bipolar,gkmw2016bending,lsw2017woods} in the so-called peanosphere topology, which is developed in \cite{matingtrees}.)

\subsection{Self-avoiding walks}
\label{subsec:saw}

\begin{figure}[ht!]
\begin{center}
\includegraphics[width=0.8\textwidth,page=2]{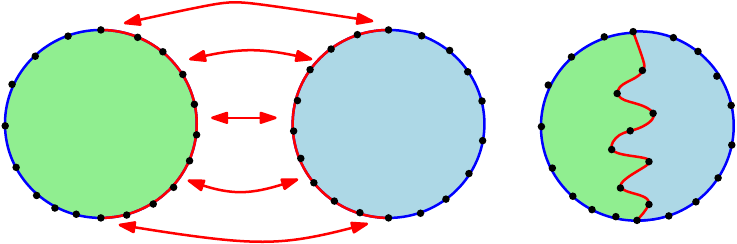}
\end{center}
\caption{\label{fig:saw} Two independent uniform quadrangulations of the disk with simple boundary, perimeter $2\ell$, and $m$ faces are glued together along a marked boundary arc of length $2s$ to produce a quadrangulation of the disk with a distinguished path of length $2s$.  The conditional law of the path is uniform among all simple paths (i.e., a self-avoiding walk) conditioned on having $m$ faces to its left and right.}
\end{figure}

Recall that the self-avoiding walk (SAW) is the uniform measure on simple paths of a given length on a graph.  The SAW on random planar maps was important historically because it was used by Duplantier and Kostov \cite{dup-kos-saw,dup-kos-saw-long} as a test case of the KPZ formula \cite{kpz1988}.  It is a particularly natural model to consider as it admits a rather simple construction.  Namely, one starts with two independent uniformly random quadrangulations of the disk with simple boundary, perimeter $2\ell$, and $m$ faces as illustrated in the left side of Figure~\ref{fig:saw}.  If one glues the two disks along a boundary segment of length $2s < 2\ell$, then one obtains a quadrangulation of the disk with perimeter $2(\ell-s)$ and $2m$ faces decorated by a distinguished path.  Conditional on the map, it is not difficult to see that the path is uniformly random (i.e., a SAW) conditioned on having $m$ faces to its left and right.  In the limit as $\ell,m \to \infty$, one obtains a gluing of two UIHPQ's with simple boundary.

The main result of \cite{gwynne-miller:uihpq} implies that the UIHPQ's converge to independent Brownian half-planes, equivalently independent $\sqrt{8/3}$-quantum wedges together with their $\QLE(8/3,0)$ metric.  Proving the convergence of the SAW in this setting amounts to showing that the discrete graph gluing of the two UIHPQ's converges to the metric gluing of the limiting Brownian half-plane instances along their boundaries.  This is accomplished in \cite{gwynne-miller:saw}.  Each Brownian half-plane instance is equivalent to a $\sqrt{8/3}$-quantum wedge with its $\QLE(8/3,0)$ metric.  In order to identify the scaling limit of the SAW with $\SLE_{8/3}$, it is necessary to show that the conformal welding of two such quantum wedges developed in \cite{she2010zipper} is equivalent to the metric gluing of the two $\sqrt{8/3}$-quantum wedges with their $\QLE(8/3,0)$ metric as it is shown in \cite{she2010zipper} that the interface is $\SLE_{8/3}$.  This is the main result of \cite{gwynne-miller:gluing}.  Combining everything gives the convergence of the SAW on random planar maps to $\SLE_{8/3}$ on $\sqrt{8/3}$-LQG.

\subsection{Percolation}
\label{subsec:perc}

\begin{figure}[ht!]
\begin{center}
\includegraphics[width=0.4\textwidth,page=5]{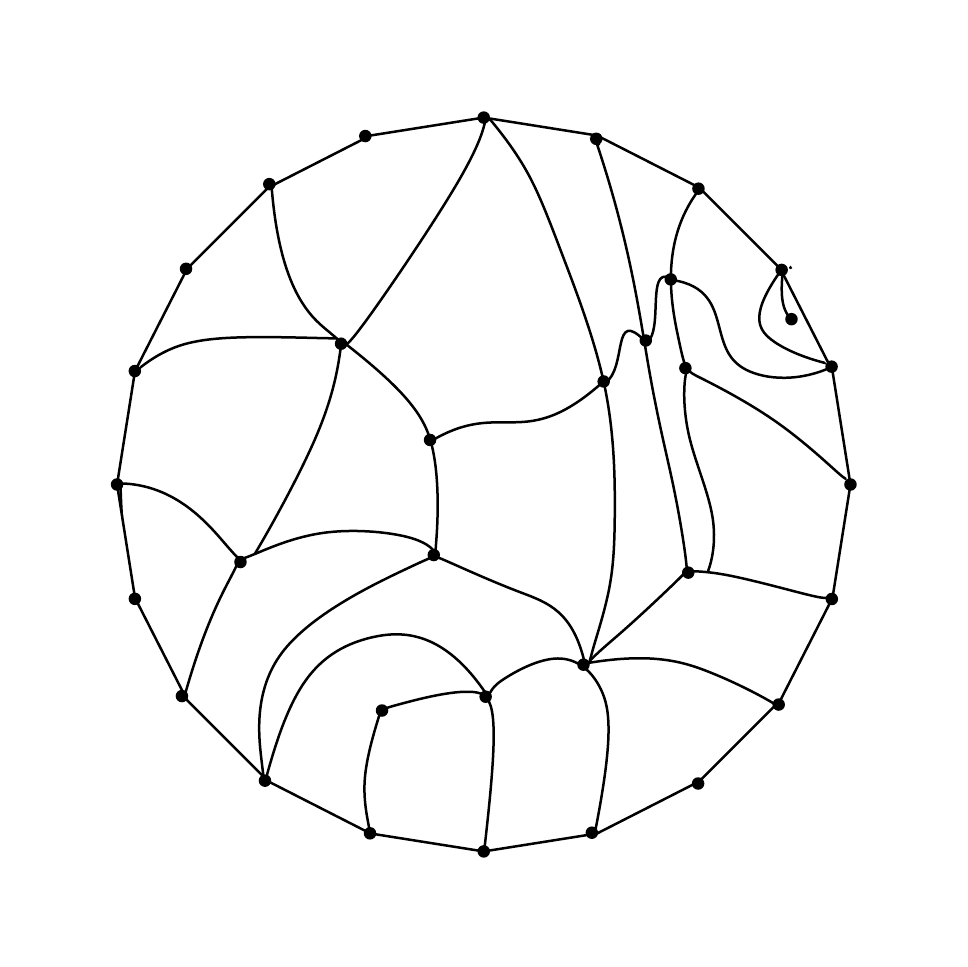}\hspace{0.03\textwidth}
\includegraphics[width=0.4\textwidth,page=4]{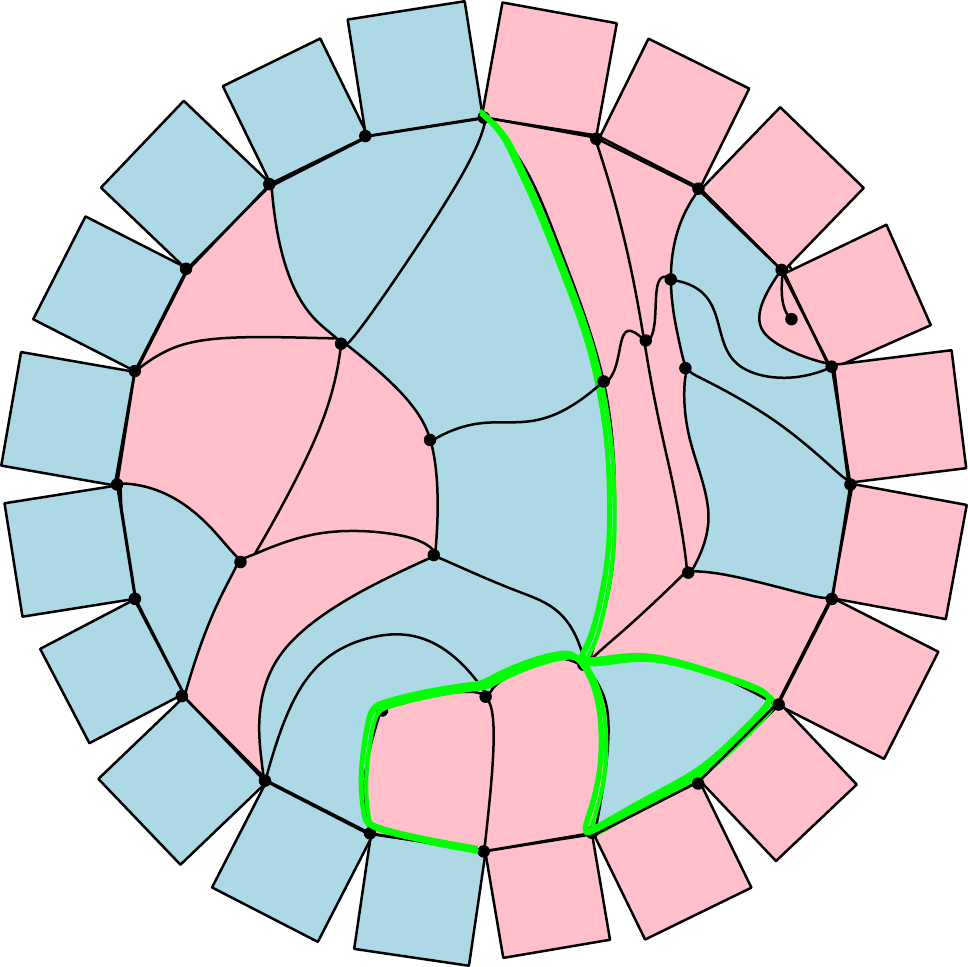}	
\end{center}
\caption{\label{fig:perc} {\bf Left:} A quadrangulation of the disk with simple boundary.  {\bf Right:} Blue (resp.\ red) quadrilaterals have been glued along two marked boundary arcs.  Faces on the inside are colored blue (resp.\ red) independently with probability $3/4$ (resp.\ $1/4$).  The green path is the interface between the cluster of blue (resp.\ red) quadrilaterals which are connected to the blue (resp.\ red) boundary arc.}
\end{figure}

It is also natural to consider critical percolation on a uniformly random planar map.  The critical percolation threshold has been computed for a number of different planar map types (see, e.g., \cite{angel-peeling,angel-curien-uihpq-perc}).  The article \cite{gm2017perc} establishes the convergence of the interfaces of critical face percolation on a uniformly random quadrangulation of the disk.  This is the model in which the faces of a quadrangulation are declared to be either open or closed independently with probability $3/4$ or $1/4$ \cite{angel-curien-uihpq-perc}.  (The reason that the critical threshold is $3/4$ and not $1/2$ is because open faces are adjacent if they share an edge and closed faces are adjacent if they share a vertex.)  The underlying quadrangulation of the disk converges to the Brownian disk \cite{bettinelli_miermont_disks,gm2017simple_disk}, equivalently a quantum disk.  The main result of \cite{gm2017perc} is that the percolation interface converges jointly with the underlying quadrangulation to $\SLE_6$ on the Brownian disk.  The proof proceeds in a very different manner than the case of the SAW.  Namely, the idea is to show that the percolation interface converges in the limit to a continuum path which is a Markovian exploration of a Brownian disk with the property that its complementary components are Brownian disks given their boundary lengths and it turns out that this property characterizes $\SLE_6$ on the Brownian disk \cite{gm2017char} (recall also Theorem~\ref{thm:sle6}).

\bibliographystyle{abbrv}
\bibliography{icm}

\end{document}